\documentclass[lettersize,journal]{IEEEtran}

\usepackage[T1]{fontenc}
\usepackage[utf8]{inputenc}
\usepackage{amsmath,amsthm,mathtools,amssymb}
\usepackage{graphicx}
\usepackage{xcolor}

\usepackage{url}

\usepackage[breaklinks]{hyperref}
\usepackage{orcidlink}
\usepackage{subcaption}

\usepackage{dutchcal}
\usepackage{dsfont}
\usepackage{siunitx}

\usepackage{tikz}
\usetikzlibrary{calc,positioning,fit}

\newcommand{\matr}[1]{{\boldsymbol{#1}}}
\renewcommand{\vec}[1]{{\boldsymbol{#1}}}

\newcommand{\eye}{\matr{I}}
\newcommand{\nn}{\mathtt{N}}
\newcommand{\imu}{{\rm j}}

\begin{document}

\title{Voltage-Sensitive Distribution Factors for Contingency Analysis and Topology Optimization}

\author{%
Maurizio Titz\orcidlink{0000-0002-0249-6244},
Dirk Witthaut\orcidlink{0000-0002-3623-5341},
Joost van Dijk\orcidlink{0009-0003-1616-3621},
Benjamin Petrick\orcidlink{0009-0001-9510-8799},
Nico Westerbeck\orcidlink{0009-0004-1761-4305},
\thanks{
M.~Titz and D.~Witthaut are with the Forschungszentrum J\"ulich, Institute of Climate and Energy Systems: Energy Systems Engineering (ICE-1),  52428 J\"ulich, Germany and the 
Institute for Theoretical Physics, University of Cologne, 50937 K\"oln, Germany. 
J.~van Dijk is with Tennet TSO BV, Utrechtseweg 310, 6812 AR Arnhem, The Netherlands.
B.~Petrick and N.~Westerbeck are with the Elia Group, Boulevard de l’Empereur 20, 1000 Brussels, Belgium.
D.~Witthaut acknowledges financial support from the Elia Group via the Elia Group Open Research Challenge 2025 and from the German Federal Ministry of Research, Technology and Space (Bundesministerium für Forschung, Technologie und Raumfahrt) via the grant number 03SF0751.
}
}

\maketitle

\begin{abstract}
Topology optimization is a promising approach for mitigating congestion and managing changing grid conditions, but it is computationally challenging and requires approximations. Conventional distribution factors like PTDFs and LODFs, based on DC power flow, fail to capture voltage variations, reactive power, and losses, thereby limiting their use  in detailed optimization tasks such as busbar splitting. 
This paper introduces generalized distribution factors derived from a voltage-sensitive linearization of the full AC power flow equations. The proposed formulation accurately reflects reactive power flows, Ohmic losses, and voltage deviations while remaining computationally efficient. We derive and evaluate generalized PTDFs, LODFs, and topology modification factors using matrix identities. We discuss potential applications including voltage-aware $N-1$ security analysis and topology optimization with a focus on busbar splitting.  Numerical experiments demonstrate close agreement with full AC solutions, significantly outperforming the traditional DC approximation.
\end{abstract}


\section{Introduction}

The transition to renewable energy sources is challenging the operation and stability of power transmission grids. Wind and solar power are fluctuating and uncertain, such that highly flexible backup generation or storage facilities are needed~\cite{entsoe_position_flex}. System operators must account for this uncertainty, for instance in contingency analysis~\cite{hasan2019existing} or voltage control \cite{knittel2020voltage}. Furthermore, renewable power sources are often built at locations with favorable environmental conditions, requiring long-distance power transmission~\cite{rodriguez2014transmission}. As a consequence, transmission system operators are facing increasing grid loads~\cite{titz2024identifying} and seek new options to mitigate congestion, for instance via dynamic topology reconfiguration~\cite{subramanian2021exploring}. Overall, system operators need new computational methods to address the challenges of the energy transition~\cite{marot2022perspectives,viebahn24}.

Distribution factors are standard tools in the design and operation of transmission grids. The standard power transfer distribution factors (PTDFs) and line outage distribution factors (LODFs) describe how line flows change when power injections are adapted or transmission lines fail~\cite{wood2013power}. Recently, distribution factors for other topology changes have been introduced for grid planning and topology optimization~\cite{goldis2016shift,van2023bus}. All these distribution factors are derived using the DC approximation, a linearized form of the power flow equations~\cite{wood2013power}, which is valid when voltages are close to their set points and losses and grid loads are moderate~\cite{purchala2005usefulness,stott2009dc}. 

In recent years, several generalizations of the DC approximation have been proposed that avoid some of the assumptions and simplifications, see Sec.~\ref{sec:literature} for a literature review. The derivation generally starts from the full nonlinear power flow equations, which include reactive power and Ohmic losses. In this article, we derive distribution factors for generalized linear power flow approximations including reactive power, losses and voltage variations. We provide expressions generalizing PTDFs, LODFs and bus merge and split distribution factors.

An additional advantage of the proposed formulation is that it retains the algebraic structure of classical distribution factors, enabling efficient low-rank updates and highly parallel implementations. In particular, the computations map naturally to GPU architectures, allowing extremely high evaluation rates for contingency screening and topology analysis.

The article is organized as follows. Sec.~\ref{sec:literature} reviews related work on generalized linear power flow approximations. In Sec.~\ref{sec:lin} we derive the linearized power flow equations and introduce a convenient vectorial notation. Sec.~\ref{sec:vlodf} develops generalized distribution factors for branch modifications and line outages. In Sec.~\ref{sec:switch} we extend the approach to switches and busbar couplers. The numerical performance is evaluated in Sec.~\ref{sec:performance} before concluding in Sec.~\ref{sec:conclusion}.

\section{Literature Review}
\label{sec:literature}

The DC approximation relies on three assumptions: (i) Ohmic losses are negligible, (ii) phase differences are small such that the sine can be linearized, (iii) voltage magnitudes are close to the set points. The validity of these approximations is studied in~\cite{purchala2005usefulness,stott2009dc,Hert06,hartmann2024synchronized}.

Several generalizations that try to avoid these assumptions have been developed:
(i) Models that incorporate Ohmic losses but keep the assumption of fixed voltage magnitudes have been proposed for instance in Refs.~\cite{stott2009dc,Simp16b}.
(ii) Methods to reduce or bound the linearization error \cite{Dorf13c,Simp16b,hartmann2024synchronized}.
(iii) Several approaches that linearize the full AC power flow equations, keeping Ohmic losses and voltage variations, have been proposed. This is the focus of the current manuscript. To arrive at a linear set of equations, most authors introduce the approximations,
\begin{align*}
    v_i v_j \sin(\theta_{i}-\theta_j) &\approx \theta_{ij}, \\
    v_i v_j \cos(\theta_{i}-\theta_j) &\approx f(v_i) - f(v_j).
\end{align*}
The function $f$ may be linear~\cite{Zhang2013}, quadratic~\cite{Yang2018}, logarithmic~\cite{Li2018}, or an arbitrary monomial~\cite{Fate15,Yang2019}. 

Distribution factors have been used for decades~\cite{landgren1972transmission}. The standard formulation is based on the DC approximation~\cite{wood2013power} and thus limited to lossless real power flow. Textbooks typically first introduce PTDFs and then treat line outages via a mapping to PTDFs~\cite{wood2013power}. This approach is highly instructive, but may lead to heavy notation and impede the generalization to more complex problems such as changes to the grid topology.

Different algebraic approaches that allow for a direct treatment of topology changes have been proposed. The DC approximation typically leads to sparse linear equations, which can be exploited in numerical computations~\cite{tinney1985sparse} and for handling topology changes~\cite{alsac1983sparsity,kim1985contingency}. Early research on sparse matrix methods was reviewed in~\cite{stott1987overview}. Other authors focus on the application of the Woodbury matrix identity, for instance in the derivation of (generalized) LODFs~\cite{guler2007generalized,guo2009direct,kaiser2020collective}. 

More recently, other topology changes have attracted research interest, as they enable optimization of the grid topology~\cite{subramanian2021exploring}. This line of research includes 
updates of phase shifting transformers~\cite{srinivasan1985psdf}, line closings~\cite{sauer2001extended}, and the merging~\cite{goldis2016shift} and splitting of busbars~\cite{van2023bus}.

\section{Linearized power flow equations}
\label{sec:lin}

\subsection{Linearizing the AC power flow equations}

For clarity, we begin by analyzing a single branch $(f,t)$. The power injections at the from end are expressed in terms of the voltage magnitudes $v$ and phase angles $\theta$ as
\begin{align}
    p_f &= v_f^2 g_{ff} + v_f v_t (g_{ft} \cos(\theta_{ft}) + b_{ft} \sin(\theta_{ft})), \nonumber \\
    q_f &= -v_f^2 b_{ff} + v_f v_t (-b_{ft}\cos(\theta_{ft}) + g_{ft}\sin(\theta_{ft})). 
    \label{eq:loadflow-singleline}
\end{align}
where we define the abbreviation $\theta_{ft}= \theta_f-\theta_t=-\theta_{tf}$. Power injections at the to end have the same form with the indices $f$ and $t$ flipped.

For the common $\pi$-line model, the admittance parameters are given by  
\begin{align*}
    \begin{pmatrix}
    g_{ff} + \imu b_{ff} & g_{ft} + \imu b_{ft} \\
    g_{tf} + \imu b_{tf} & g_{tt} + \imu b_{tt}
    \end{pmatrix}
    = 
    \begin{pmatrix}
    \frac{y_{ft}^s + y^{cf}_{ft}}{\tau^2}  &
    \frac{- y_{ft}^s}{\tau \exp(-\imu \alpha_{ft})} \\
    \frac{- y_{ft}^s}{\tau \exp(+\imu \alpha_{ft})} &
    y_{ft}^s + y^{ct}_{ft} 
    \end{pmatrix},
\end{align*}
where $\tau$ denotes the tap ratio, $y_{ft}^s$ denotes the series admittance, and the parameters $y^{cf}_{ft}$ and $y^{ct}_{ft}$ account for the charging elements. Furthermore, $\alpha_{ft}= - \alpha_{tf}$ denotes the phase angle of a phase shifting transformer; for an ordinary line or transformer $\alpha_{ft}=0$.

We propose a linearization that is particularly useful for the analysis of line outages and topology modifications. Assume that the voltage phase angle difference is close to a reference value $\hat \theta_{ft}$. Then we can linearize the trigonometric functions,
\begin{align*}
    \cos(\theta_{ft}) &\approx \cos(\hat \theta_{ft}) 
    - \sin(\hat \theta_{ft}) 
       (\theta_{ft} - \hat \theta_{ft}) , \\
    \sin(\theta_{ft}) & \approx \sin(\hat \theta_{ft}) 
    + \cos(\hat \theta_{ft}) (\theta_{ft} - \hat \theta_{ft}).
\end{align*}
We can follow two strategies for the definition of the reference state: (i) For general computations, we can use a flat initial guess where $\hat \theta_{ft}=0$ for ordinary transmission lines and  $\hat \theta_{ft}=\alpha_{ft}$ for PSTs. This strategy is commonly referred to as a cold-start model. (ii) To assess line outages and topology
modifications, we can use the pre-modification state for $\hat \theta_{ft}$. This strategy is commonly referred to as a hot-start model.

For each bus $i$, the voltage magnitude is written as
$
    v_i = \hat v_i (1 + u_i)
$
where $u_i$ is a dimensionless small parameter. For slack and PV nodes, we use the reference values for the voltage magnitudes $\hat v_i = v_{i,\rm set}$. For the PQ nodes we use (i) 
the reference value of the respective voltage level in a cold-start model or (ii) the pre-modification state in a hot-start model.

We can then substitute this approach into Eq.~\eqref{eq:loadflow-singleline} and neglect all nonlinear terms to obtain a linear relation between real and reactive power and the state variables $\theta$ and $u$. We provide these relations in a convenient matrix form which will facilitate the formulation of the network equations later. We obtain
\begin{align}
    \begin{pmatrix}
        p_f \\ p_t \\ q_f \\ q_t
    \end{pmatrix}
    = 
    \begin{pmatrix}
        \hat p_f \\ \hat p_t \\ \hat q_f \\ \hat q_t
    \end{pmatrix}
    +
    \begin{pmatrix}
    M^{p \theta}_{ff} & M^{p \theta}_{ft} & 
            M^{p u}_{ff} & M^{p u}_{ft} \\
    M^{p \theta}_{tf} & M^{p \theta}_{tt} & 
            M^{p u}_{tf} & M^{p u}_{tt} \\
    M^{q \theta}_{ff} & M^{q \theta}_{ft} & 
            M^{q u}_{ff} & M^{q u}_{ft} \\
    M^{q \theta}_{tf} & M^{q \theta}_{tt} & 
            M^{q u}_{tf} & M^{q u}_{tt} \\        
    \end{pmatrix}
    \begin{pmatrix}
        \theta_f  \\ \theta_t  \\ 
          u_f \\ u_t
    \end{pmatrix}
    \label{eq:linearization-ft}
\end{align}
with the elements
\begin{align*}
    \hat p_f &= 
       \hat{v}_f^2 g_{ff} + \hat{v}_f\hat{v}_t 
       \left(g_{ft}\cos(\hat \theta_{ft}) + b_{ft} \sin(\hat \theta_{ft}) \right) \\
       & \quad - \hat{v}_f\hat{v}_t  \left(-g_{ft}\sin(\hat \theta_{ft}) + b_{ft} \cos(\hat \theta_{ft}) \right) \hat \theta_{ft} \\
    \hat q_f &=  
        - \hat{v}_f^2 b_{ff} + \hat{v}_f\hat{v}_t
        \left( g_{ft}  \sin(\hat \theta_{ft})  -b_{ft} \cos(\hat \theta_{ft}) \right) \\
      & \quad -   \hat{v}_f\hat{v}_t  \left(b_{ft}\sin(\hat \theta_{ft}) + g_{ft} \cos(\hat \theta_{ft}) \right) \hat \theta_{ft} \\
    M^{pu}_{ff} &= 2 \hat{v}_f^2 g_{ff} 
        + \hat{v}_f\hat{v}_t
        \left(g_{ft}\cos(\hat \theta_{ft}) + b_{ft} \sin(\hat \theta_{ft}) \right)  \\
    M^{pu}_{ft} &= 
        + \hat{v}_f\hat{v}_t
        \left(g_{ft}\cos(\hat \theta_{ft}) + b_{ft} \sin(\hat \theta_{ft}) \right)  \\
    M^{qu}_{ff} &=     
        - 2 \hat{v}_f^2 b_{ff} 
        + \hat{v}_f\hat{v}_t
        \left( -b_{ft} \cos(\hat \theta_{ft}) + g_{ft}          \sin(\hat \theta_{ft}) \right) \\
    M^{qu}_{ft} &=     
        + \hat{v}_f\hat{v}_t
        \left( -b_{ft} \cos(\hat \theta_{ft}) + g_{ft}          \sin(\hat \theta_{ft}) \right) \\   
    M_{ff}^{p\theta} &= - M_{ft}^{p\theta}
        = \hat{v}_f\hat{v}_t \left( -g_{ft}\sin(\hat \theta_{ft}) + b_{ft}\cos(\hat \theta_{ft}) \right) \\
    M_{ff}^{q\theta} &= - M_{ft}^{q\theta}
    = \hat{v}_f\hat{v}_t
    \left( b_{ft}\sin(\hat \theta_{ft}) + g_{ft}\cos(\hat \theta_{ft}) 
    \right) .
\end{align*}
The remaining matrix elements for the to end of the branch have the same form, with all indices $f$ and $t$ swapped.

\subsection{Shunt elements}

Shunt elements are described by an admittance $y^{\text{shunt}}_f = g^{\text{shunt}}_f + \imu b^{\text{shunt}}_f$ connecting a bus $f$ to the ground. Using the linearization introduced above, they contribute to the power balance according to the equation
\begin{align}
    \begin{pmatrix}
        p_f \\ q_f 
    \end{pmatrix}
    = 
    \begin{pmatrix}
        {\hat p}_f^{\text{shunt}} \\ 
        \hat q_f^{\text{shunt}}
    \end{pmatrix}
    +
    \begin{pmatrix}
        N^{p \theta}_{f} & N^{p u}_{f} \\
        N^{q \theta}_{f}    & N^{q u}_{f} \\       
    \end{pmatrix}
    \begin{pmatrix}
        \theta_f \\  u_f
    \end{pmatrix}
     \label{eq:linearization-shunt}
\end{align}
with the parameters
\begin{align*}
    \hat p_f^{\text{shunt}} &= \hat v_f^2 g^{\text{shunt}}_f, &
    N^{p u}_{f} &= 2 \hat v_f^2 g^{\text{shunt}}_f,
    \\
    \hat q_f^{\text{shunt}} &= -\hat v_f^2 b^{\text{shunt}}_f, &
    N^{q u}_{f} &= - 2 \hat v_f^2 b^{\text{shunt}}_f 
\end{align*}
and $N^{p \theta}_{f} = N^{q \theta}_{f} =0$.

\subsection{The network equations}

We now collect all prior results to formulate the linearized power flow equations for an extended grid. We number the PQ nodes as $i = 1,\ldots,n$, PV nodes as $i = n+1,\ldots,n+m$ and the slack node as $n+m+1$. The total number of nodes is denoted as $\nn = n+m+1$.

We introduce a convenient vectorial notation, summarizing all power injections and state variables in the vectors
\begin{align*}
    \vec p &= ( p_1  , p_2, \ldots,   p_{n+m} )^\top,
    &
    \vec q &= ( q_1  , q_2, \ldots,   q_n  )^\top
    \\
    \vec \theta &= (\theta_1,\theta_2, \ldots, \theta_{n+m} )^\top,
    & 
    \vec u &= (u_1, u_2 , \ldots,  u_n )^\top .
\end{align*}

For bookkeeping we define the index vectors $\vec \eta, \vec \zeta \in \mathbb{R}^{2n+m}$ as follows. The vector $\vec \eta_i$ has a single unit entry at position $i$ and the vector $\vec \zeta_i$ has a single unit entry at position $i+n+m$. All other entries vanish. We note that the vector $\vec \eta_i$ is linked to the real power or phase angle of a bus $i$, whereas the vector $\vec \zeta_i$ is linked to the reactive power or voltage magnitude of a bus $i$. Since we do not have to include these quantities for PV nodes, we set $\vec \zeta_i = \vec 0$ for $i>n$. Furthermore, we define the shorthands $\vec \mu_{ij} = \vec \eta_i - \vec \eta_j$ and $\vec \nu_{ij} = \vec \zeta_i - \vec \zeta_j$.

The parameters introduced in the linearization in Eq.~\eqref{eq:linearization-ft} and Eq.~\eqref{eq:linearization-shunt} are summarized in two vectors with constant elements
\begin{align}
    \vec{\hat p} &= 
       \sum_{f=1}^{n+m} \left( \hat p_f + \hat p_f^{\rm shunt} \right)  \vec \eta_f  
    &
    \vec{\hat q} &= 
       \sum_{f=1}^n \left( \hat q_f + \hat q_f^{\rm shunt} \right)  \vec \zeta_f     .
    \label{eq:pqhat-fullvector}
\end{align}
Then we collect all matrix elements describing branches and shunt elements in the matrix 
\begin{align}
    \label{eq:MN-fullmatrix}
    &\matr M = \sum_{f=1}^{m+n}
    N^{p u}_{f}   \vec \eta_f \vec \zeta_f^\top
    +N^{q u}_{f} \vec \zeta_f \vec \zeta_f^\top.
    \\
    & +
    \sum_{f,t=1}^{m+n} 
    M^{p\theta}_{f t} \vec \eta_f \vec \eta_t^\top
    +M^{p u}_{f t} \vec \eta_f \vec \zeta_t^\top
    +M^{q\theta}_{f t} \vec \zeta_f \vec \eta_t^\top
    +M^{q u}_{f t} \vec \zeta_f \vec \zeta_t^\top .
    \nonumber
\end{align}
We note that the matrix $\matr M$ corresponds to the Jacobian of the AC power flow equations at the reference state. We thus obtain the linearized power flow equations in a compact vectorial form
\begin{align}
    \begin{pmatrix}
       \vec p \\ \vec q 
    \end{pmatrix}
    =
    \begin{pmatrix}
       \vec{\hat p} \\ 
       \vec{\hat q} 
    \end{pmatrix}
    +
    \matr M
    \begin{pmatrix}
       \vec \theta \\ \vec u
    \end{pmatrix}.
    \label{eq:loadflow-lin-vector}
\end{align}
We will refer to this as the generalized or voltage-sensitive DC+ approximation in the following.

\section{Finite branch modifications and line outages}
\label{sec:vlodf}

In this section, we consider arbitrary finite branch modifications which may include line outages, line closings or changes of the PST phase angle. Switches as well as bus split and merge operations will be discussed in Sec.\ref{sec:switch}.

\subsection{Branch modification distribution factors}
\label{sec:branch-finite-general}

We start from a reference grid for which the linearized power flow equations are given by Eq.~\eqref{eq:loadflow-lin-vector}. Notably, if we have solved the full nonlinear power flow equations for this reference grid, we may use the resulting voltage magnitudes as the reference values $\hat v_i$ and $\hat \theta_i$.

We consider the modification of a single branch $(f,t)$, leading to changes to the matrix $\matr M$ and the inhomogeneities,  
\begin{align*}
    \matr M &\rightarrow {\matr M}' = \matr M + \Delta \matr M \\
    \vec {\hat p} &\rightarrow \vec{\hat p}' =  
    \vec{\hat p} + \Delta \vec{\hat p},
    &
    \vec{\hat q} & \rightarrow \vec{\hat q}' =  
      \vec{\hat q} + \Delta \vec{\hat q} .
\end{align*}
The matrix $\Delta \matr M$ is of rank three or less and is written as 
\begin{align}
    \Delta \matr M 
    &=  \underbrace{
    \begin{pmatrix}
       \vec{\mathcal{s}}_{ft} & \vec s_f & \vec s_t
    \end{pmatrix}
    }_{=: \matr S}
    \underbrace{
    \begin{pmatrix}
        \vec \eta_f^\top - \vec \eta_t^\top \\
        \vec \zeta_f^\top \\
        \vec \zeta_t^\top
    \end{pmatrix}
    }_{=: \matr R}
    \label{eq:Delta-M-single}
\end{align}
with
\begin{align*}
    \vec s_f = \Delta M^{pu}_{f f}  \vec \eta_f +
           \Delta M^{pu}_{t f}  \vec \eta_t + 
           \Delta M^{qu}_{f f}  \vec \zeta_f +
           \Delta M^{qu}_{t f}  \vec \zeta_t \\
    \vec s_t = \Delta M^{pu}_{f t}  \vec \eta_f +
           \Delta M^{pu}_{t t}  \vec \eta_t + 
           \Delta M^{qu}_{f t}  \vec \zeta_f +
           \Delta M^{qu}_{t t}  \vec \zeta_t \\
    \vec{\mathcal{s}}_{ft} = 
            \Delta M^{p\theta}_{f f} \vec \eta_f +
            \Delta M^{p\theta}_{t f} \vec \eta_t +
            \Delta M^{q\theta}_{f f} \vec \zeta_f +
            \Delta M^{q\theta}_{t f} \vec \zeta_t .
\end{align*}
Here, we have used the fact that the shunt elements are not affected and that $M^{p\theta}_{ff} = -M^{p\theta}_{tf}$ and $M^{q\theta}_{ff} = - M^{q\theta}_{tf}$ and similar for $f$ and $t$ swapped. The change in topology induces a change of the state variables 
\begin{align*}
    \vec \theta \rightarrow {\vec \theta}' &= \vec \theta + \vec \psi, 
    &
    \vec u \rightarrow {\vec u}' &= \vec u + \vec \sigma.
\end{align*}
To compute $\vec \psi$ and $\vec \sigma$, we subtract the linear power flow equations \eqref{eq:loadflow-lin-vector} before and after the  modification and obtain
\begin{align}
    \begin{pmatrix}
       \Delta \vec{\hat p} \\ \Delta \vec{\hat q} \end{pmatrix}
    +
    \Delta \matr M
    \begin{pmatrix}
       \vec \theta \\ \vec u 
    \end{pmatrix}
    = 
    - \left( \matr M +  \Delta \matr M \right)
    \begin{pmatrix}
       \vec \psi \\ \vec \sigma
    \end{pmatrix},
\end{align}
which is formally solved as
\begin{align}
    \begin{pmatrix}
       \vec \psi \\ \vec \sigma
    \end{pmatrix}
    =
    - \left( \matr M +  \Delta \matr M \right)^{-1}
    \left[ 
    \begin{pmatrix}
       \Delta \vec{\hat p} \\ \Delta \vec{\hat q} 
    \end{pmatrix}
    +
    \Delta \matr M
    \begin{pmatrix}
       \vec \theta \\ \vec u 
    \end{pmatrix}
    \right].
    \label{eq:psi-sigma}
\end{align}
Since $\Delta \matr M$ has low rank, we can compute the inverse of $\matr M + \Delta \matr M$ efficiently using the Woodbury matrix identity,
\begin{align}
    &\left( \matr M + \Delta \matr M \right)^{-1} 
    \nonumber \\
    &= \matr M^{-1} - 
    \matr M^{-1}  \matr S 
    \left( \eye + \matr R \matr M^{-1} \matr S \right)^{-1}
    \matr R \matr M^{-1} ,
    \label{eq:woodbury}
\end{align}
where $\eye$ denotes the identity matrix.
We emphasize that the matrix $\left( \eye + \matr R \matr M^{-1} \matr S \right)$ has size $3 \times 3$, such that we can easily compute its inverse. This approach is particularly important if we want to evaluate various branch modifications, as for instance all transmission line outages for a comprehensive security assessment. In this case we have to invert $\matr M$ only once for the reference grid topology. For any branch modification we can then compute $\left( \matr M + \Delta \matr M \right)^{-1} $ efficiently using the Woodbury matrix identity.

\subsection{Multiple branch modifications}

The above approach is easily generalized to multiple branch modifications or multiple branch outages. We label the modified branches as $(f_1,t_1), (f_1,t_2), \ldots, (f_k,t_k)$. The change of the matrix $\matr M$ can then be written as
\begin{align*}
    \Delta \matr M &=  
    \sum_{i=1}^k \begin{pmatrix}
       \vec{\mathcal{s}}_{f_i t_i} & \vec s_{f_i} & \vec s_{t_i}
    \end{pmatrix}
    \begin{pmatrix}
        \vec \eta_{f_i}^\top - \vec \eta_{t_i}^\top \\
        \vec \zeta_{f_i}^\top \\
        \vec \zeta_{t_i}^\top
    \end{pmatrix}
    \nonumber \\
    &= \matr S \matr R
\end{align*}
and is of rank $3k$. The change in the state variables is still given by Eq.~\eqref{eq:psi-sigma}, and the inverse can be computed using the Woodbury matrix identity in Eq.~\eqref{eq:woodbury}, where the matrix $\left( \eye + \matr R \matr M^{-1} \matr S \right)$ now has size $3k \times 3k$.

\subsection{Line modification and outage distribution factors}

We now consider the important special case of an ordinary transmission line $(f,t)$ that satisfies
\begin{align*}
    \hat \theta_{ft} \approx 0, \qquad
    \hat v_{f} \approx \hat v_t,  \qquad
    \alpha_{ft}=0, \qquad
    \tau = 1.
\end{align*}
We consider a modification of the series admittance,
\begin{align}
    b^s_{ft} \rightarrow  b^s_{ft} + \Delta b^s_{ft}, \qquad
    g^s_{ft} \rightarrow  g^s_{ft} + \Delta g^s_{ft},
\end{align}
while we assume that charging elements are not modified. Notably, a line outage is described by $\Delta b^s_{ft}= - b^s_{ft}$ and $\Delta g^s_{ft}= - g^s_{ft}$.
In this case the matrix $\Delta \matr M$ defined in Eq.~\eqref{eq:Delta-M-single} simplifies considerably. It is of rank 2 and can be written as
\begin{align*}
    \Delta \matr M = 
    \begin{pmatrix}
        \vec \mu_{ft}  &
        \vec \nu_{ft}
    \end{pmatrix}
    \underbrace{\hat v_f \hat v_t
    \begin{pmatrix}
        - \Delta b^s_{ft} & 
        + \Delta g^s_{ft} \\
        - \Delta g^s_{ft} &
        - \Delta b^s_{ft}
    \end{pmatrix}
    }_{=: \matr D}
    \begin{pmatrix}
        \vec \mu_{ft}^\top \\
        \vec \nu_{ft}^\top
    \end{pmatrix},
\end{align*}
using the short-hands $\vec \mu_{ft} =  \vec \eta_f - \vec \eta_t$ and $\vec \nu_{ft} = \vec \zeta_f - \vec \zeta_t$. The inhomogeneities do not change so that $\Delta \vec{\hat p} = \Delta \vec{\hat q} = 0$. 

The inverse $(\matr M + \Delta \matr M)^{-1}$ can now be computed from $\matr M^{-1}$ using Eq.~\eqref{eq:woodbury}, requiring the inversion of only a $2\times2$ matrix. The changes of the state variables $\vec \psi$ and $\vec \sigma$ are then readily obtained using Eq.~\eqref{eq:psi-sigma}. If we monitor another branch $(k, l)$ in the grid, we can project the resulting equations onto $\vec \mu_{kl}$ or $\vec \nu_{kl}$ to obtain
\begin{align}
    \begin{pmatrix}
        \psi_k - \psi_l  \\  \sigma_k - \sigma_l
    \end{pmatrix} =
    LMDF_{(kl),(ft)}     
    \begin{pmatrix}
        \theta_f - \theta_t    \\   u_f - u_t
    \end{pmatrix} 
    \label{eq:Delta-from-VLODF}
\end{align}
with the voltage-sensitive line modification distribution factor
\begin{align}
    \label{eq:VLODF}
     LMDF_{(kl),(ij)} &= 
    \matr A \left( \eye + \matr D \matr A \right)^{-1} \matr D
    \\
    \matr A &= 
    \begin{pmatrix}
        \vec \mu_{kl}^\top \matr M^{-1} \vec \mu_{ft} &
        \vec \mu_{kl}^\top \matr M^{-1} \vec \nu_{ft} \\
        \vec \nu_{kl}^\top \matr M^{-1} \vec \mu_{ft} &
        \vec \nu_{kl}^\top \matr M^{-1} \vec \nu_{ft} 
    \end{pmatrix}. 
\end{align}
This result can be compared to the standard line modification or line outage distribution factors (LODFs) used extensively in power system stability analysis. In both cases, the impact of a line modification is essentially described by the inverse of a matrix defined by the grid topology -- the nodal susceptance matrix in the case of ordinary LODFs and the matrix $\matr M$ in our case. The main difference is that the generalized distribution factor in Eq.~\eqref{eq:VLODF} is not a number, but a $2\times2$ matrix. This is necessary to describe the coupling between the two state variables, the phase angles and voltage magnitudes.

\subsection{Numerical evaluation}

We test the accuracy of the generalized DC+ approximation by comparing it to full AC simulations. We use plain AC power flow computations in PyPowSyBl~\cite{powsybl2025} without any adaptable grid elements. That is, we disable any type of controller, generator limits, distributed slack and DC grid elements and ensure that the $\pi$-line model is used for all transformers.

In each test we first carry out an AC power flow for the $N-0$ case and use the results for the reference values $\hat v_i$ and $\hat \theta_{ft}$ and the computation of the matrix $\matr M$. Then we run an $N-1$ analysis, simulating the failure of each branch in the grid. For each failure, we compute the state variables $v_m$, $\theta_m$ and the power injections $P_m$ and $Q_m$ for each bus $i$ via the generalized distribution factors given by Eq.~\eqref{eq:psi-sigma}. We then evaluate the approximation errors, i.e. the deviation from the full AC results. In addition, we compute the approximation error for the common DC approximation for comparison.

\begin{figure}[tb]
    \centering
    \includegraphics[width=\columnwidth]{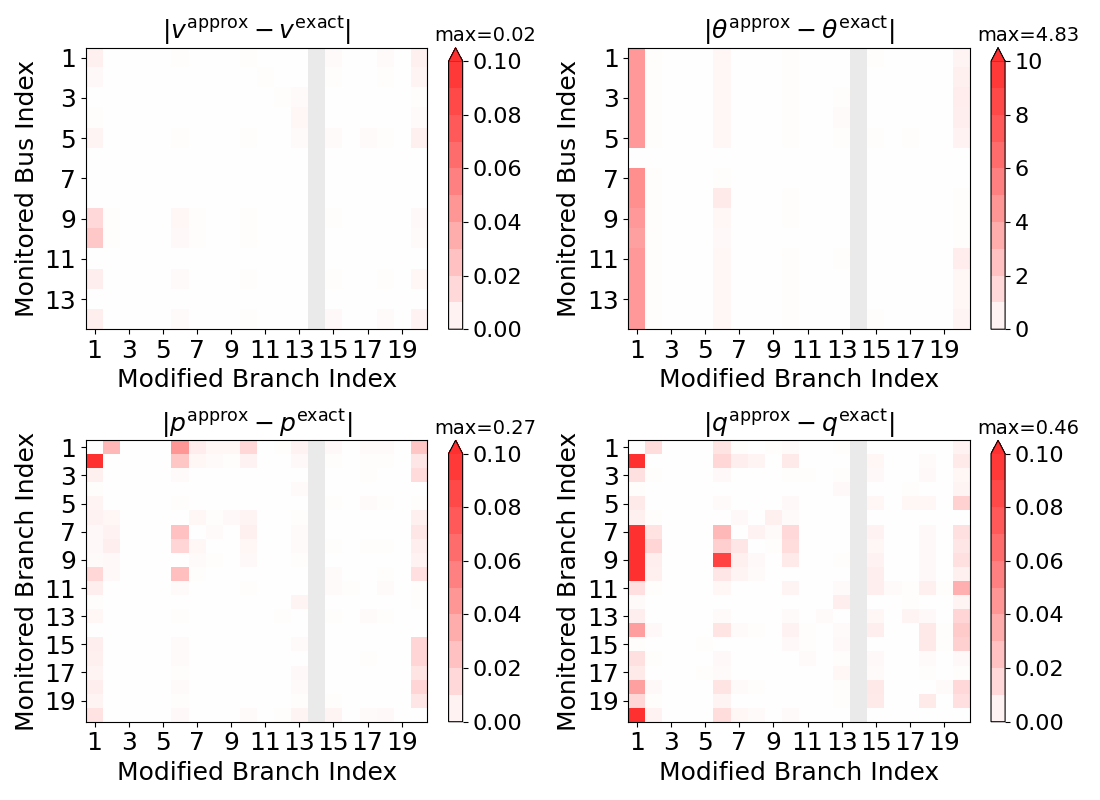}
    \caption{Approximation errors of the generalized DC+ approximation defined by Eq.~\eqref{eq:psi-sigma} with respect to the full AC simulation for the IEEE 14-bus test grid. AC simulations and grid data from PyPowSyBl. Branch 14 is excluded since its outage would disconnect the grid.
    }
    \label{fig:error-case14}
\end{figure}

A detailed analysis for the IEEE 14-bus test grid is provided in Fig.~\ref{fig:error-case14}. We find that the approximation error in the state variables is very small for almost all outage scenarios; with moderate errors for the outage of branch 1. Most importantly, the generalized DC+ approximation given by Eq.~\eqref{eq:psi-sigma} captures the changes in voltage magnitudes well. Errors are most pronounced for nodal reactive power injections.

\begin{figure}[tb]
    \centering
    \includegraphics[width=\columnwidth]{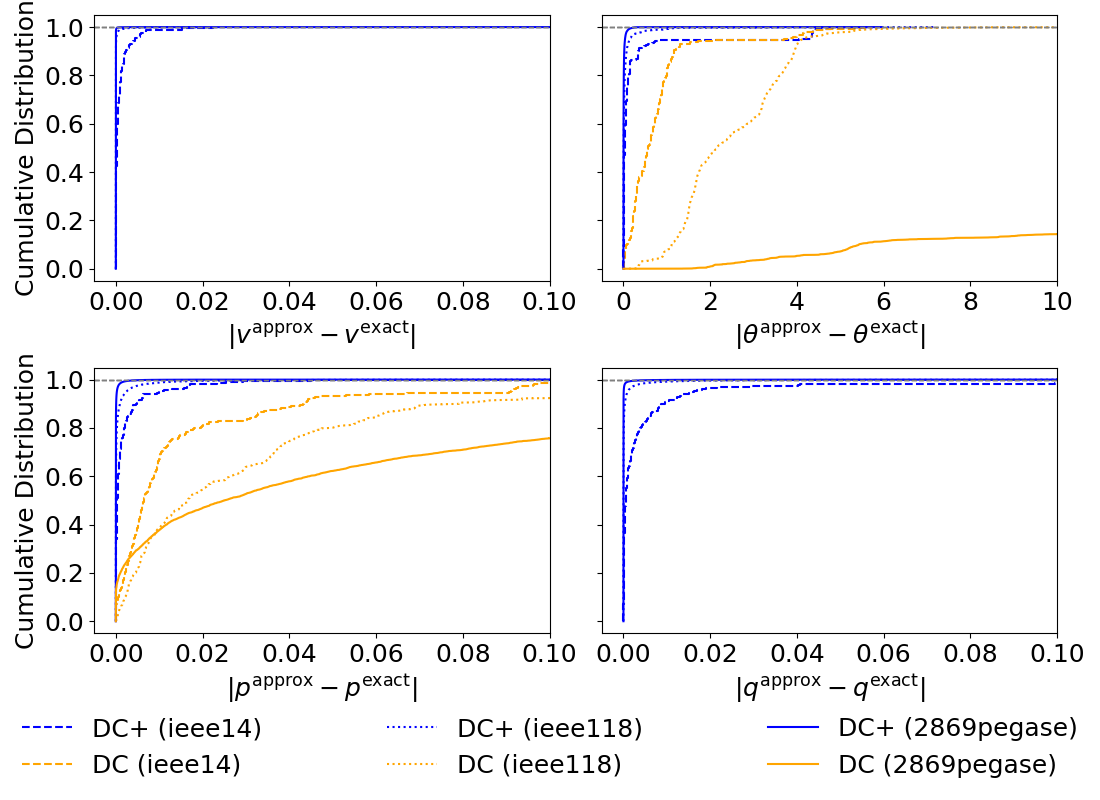}
    \caption{Statistical analysis of the approximation error for three common test grids of different sizes considering all branch outages and all nodes.
    The approximation errors of the generalized DC+ approximation defined by Eq.~\eqref{eq:psi-sigma} (blue lines) are substantially smaller than those of the standard DC approximation (orange line).
    }
     \label{fig:error-CDF}
\end{figure}

A statistical analysis of the approximation error is provided in Fig.~\ref{fig:error-CDF}. We compute the approximation errors for all branch outages and all nodes and display the cumulative distribution function. 
We observe that approximation errors of the generalized DC+ approximation~\eqref{eq:psi-sigma} are generally very small.

For the angles $\theta_m$ and real power injections $P_m$, we also show the error of the standard DC approximation for reference, computed using PowSyBl's DC security analysis~\cite{powsybl2025}. We find that the generalized DC+ approximation vastly outperforms the standard DC approximation in terms of accuracy, while the computational cost is only slightly higher.
Notably, the difference is largest for the Pegase2869 test grid, where the standard DC approximation yields large errors. However, using a distributed slack would improve the standard DC approximation such that the difference would be smaller.

\section{Switches and busbar couplers}
\label{sec:switch}

Dynamic topology reconfiguration is a promising option for reducing grid congestion~\cite{subramanian2021exploring}. This includes, in particular, opening or closing switches or busbar couplers. In the following, we will assume that the inverse of the matrix, $\matr M^{-1}$, has been computed once for a reference topology.
We show how to efficiently update this matrix after a switch or busbar coupler is modified.  

We consider an ideal switch or busbar coupler which is modeled as a branch $e = (f,t)$ with two possible states:
\begin{align}
    b_{e}^s = 0 \; \text{(open)}
    \qquad
    b_{e}^s \rightarrow \infty \; \text{(closed)}
\end{align}
while we have $g_{e}^s = \alpha_{e} = y_{e}^c = 0$ and $\tau=1$ in both cases. Furthermore, we assume that $\hat v_f=\hat v_t$. For brevity, we use $e$ as a shorthand for the branch $(f,t)$.
The limit of infinite series susceptance $b_{e}^s$ makes the treatment challenging. Even for the established DC approximation, bus split distribution factors were introduced only recently \cite{van2023bus,van2024unified}.

\subsection{Closing of switches and busbar couplers}
\label{sec:bus-merge}

We base our analysis on the results of
Sec.~\ref{sec:branch-finite-general}. We first assume that $b_e^s$ is finite and take the limit $b_e^s \rightarrow \infty$ only in the final result. In the following, the symbol $\matr M_o$ denotes the grid with the switch open and $\matr M_c$ the grid with the switch closed. We thus have $\matr M_c = \matr M_o + \Delta \matr M$ with 
\begin{align}
    \Delta \matr M 
    &=  b_e^s \hat v_t \hat v_f \;
    \underbrace{
    \begin{pmatrix}
       \vec{\mathcal{s}}_{ft} & \vec s_f & \vec s_t
    \end{pmatrix}
    }_{=: \matr S}
    \underbrace{
    \begin{pmatrix}
        \vec \eta_f^\top - \vec \eta_t^\top \\
        \vec \zeta_f^\top \\
        \vec \zeta_t^\top
    \end{pmatrix}
    }_{=: \matr R}
     \label{eq:DeltaM-switch-closing}
\end{align}
and
\begin{align*}
    \vec s_f &= 
       + \sin(\hat \theta_e) (\vec \eta_t - \vec \eta_f)
       + \cos(\hat \theta_e) (\vec \zeta_t + \vec \zeta_f)
       - 2 \vec \zeta_f \\
    \vec s_t &= 
       + \sin(\hat \theta_e) (\vec \eta_t - \vec \eta_f)
       + \cos(\hat \theta_e) (\vec \zeta_t + \vec \zeta_f)
       - 2 \vec \zeta_t \\
    \vec{\mathcal{s}}_{ft} &= 
        - \sin(\hat \theta_e) (\vec \zeta_f - \vec \zeta_t)
        - \cos(\hat \theta_e) (\vec \eta_f - \vec \eta_t).
\end{align*}

Using again the Woodbury matrix identity we find that the inverse $\matr M_c^{-1}$ can be written as
\begin{align*}
    \matr M_c^{-1} =& 
    \left[ \matr M_o + b_e^s \hat v_f \hat v_t \matr S \matr R \right]^{-1}  
    \\
    =& \matr M_o^{-1} - \matr M_o^{-1} \matr S
    \left[ 
        \frac{1}{b_e^s \hat v_f \hat v_t} \eye + 
        \matr R \matr M_o^{-1} \matr S
    \right]^{-1} 
    \matr R \matr M_o^{-1}
\end{align*}
At this point we can carry out the limit $b^s_e \rightarrow \infty$ and obtain
\begin{align}
    \matr M_c^{-1} =& 
    \matr M_o^{-1} - \matr M_o^{-1} \matr S
    \left[ 
         \matr R \matr M_o^{-1} \matr S
    \right]^{-1} 
    \matr R \matr M_o^{-1} \, .
    \label{eq:Minv-merge}
\end{align}
We can now compute the state variables for the grid with the closed switch by solving Eq.~\eqref{eq:loadflow-lin-vector} and carry out an $N-1$ contingency analysis via Eq.~\eqref{eq:psi-sigma} using $\matr M_c^{-1}$ instead of $\matr M_o^{-1}$.

\subsection{Bus splits}

Bus splits are typically more important than bus mergers because busbar couplers are typically closed in the reference grid layout.
In the following, we introduce a formulation of bus splits for a single busbar coupler. Initially, the coupler is closed such that the two busbars effectively constitute a single node of the network. We want to describe the impact of the opening of the coupler making the two busbars two separate nodes with no direct connection. For simplicity we treat the two terminal busbars as PQ nodes without shunt elements.

The treatment of a bus split involves several technical difficulties. First, the limit $b_e^s \rightarrow \infty$ is significantly more involved than in the case of a bus merge. Second, two coupled busbars effectively form a single node of the network. Starting from this single-node description, we have to add a new node to describe the opening of the busbar. 

\begin{figure}
    \centering
    \begin{subfigure}[b]{0.15\textwidth}
        \centering
        \includegraphics[width=\textwidth]{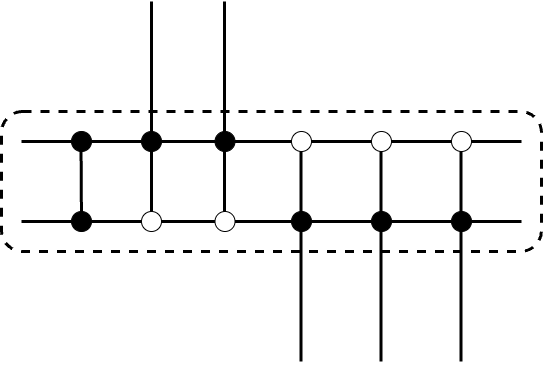}
        \caption{Merged $M_m$}
    \end{subfigure}
    \hfill
    \begin{subfigure}[b]{0.15\textwidth}
        \centering
        \includegraphics[width=\textwidth]{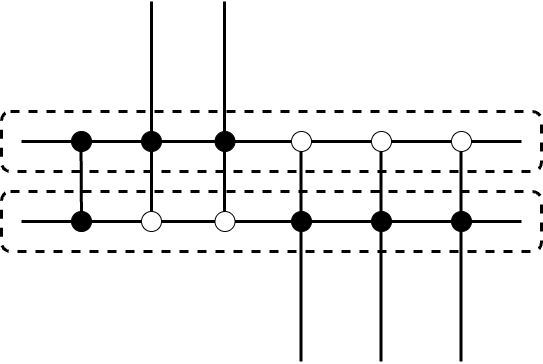}
        \caption{Closed $M_c$}
    \end{subfigure}
    \hfill
    \begin{subfigure}[b]{0.15\textwidth}
        \centering
        \includegraphics[width=\textwidth]{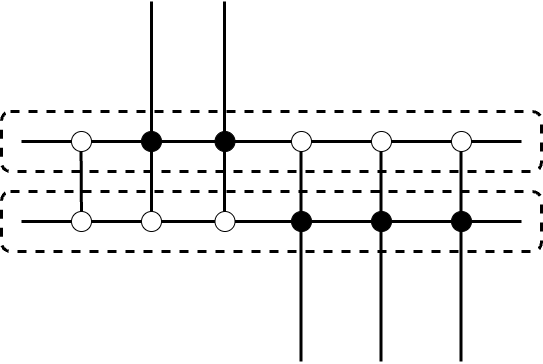}
        \caption{Open $M_o$}
    \end{subfigure}
\caption{The different grid configurations used in the treatment of the opening of a switch or busbar coupler. }
\label{fig:busbar-sketch}
\end{figure}

Let the matrix $\matr M_m$ describe the reference grid where all switches or busbar couplers are closed and merged into  single nodes. Furthermore, let $\matr M_o$ describe the grid where the busbar coupler has been opened. The two uncoupled busbars then constitute separate nodes so that the dimension of $\matr M_o$ is larger than the dimension of $\matr M_m$. To relate the case of interest to the reference case, we introduce a third network configuration $\matr M_c$ where all switches and busbar couplers are closed, but where the nodes have \emph{not} been merged. In this configuration, the incident branches must be assigned to the node (busbar), to which they are connected when the coupler opens. The two matrices $\matr M_c$ and $\matr M_o$ have the same dimension. These configurations are sketched in Fig.~\ref{fig:busbar-sketch}.

Our computation is based on the inverse matrices $\matr M^{-1}_{m,c,o}$, assuming that $\matr M^{-1}_{m}$ is computed once for the reference grid. We proceed via the intermediate grid configuration. The inverse matrix $\matr M^{-1}_{c}$ is obtained from $\matr M^{-1}_{m}$ by a simple padding procedure. That is, the rows and columns of the matrix $\matr M^{-1}_{m}$ corresponding to the coupled busbars are simply copied and inserted at the correct position. The only practical challenge for implementations is correct bookkeeping. 
At this step, the vectors $\vec p$ and $\vec q$ as well as $\vec{\hat p}$ and $\vec{\hat q}$ must be updated in a consistent way. Power injections must be assigned to the node to which they are connected after the split.

We now show how to compute $\matr M^{-1}_{o}$ from $\matr M^{-1}_{c}$. As before, we first assume $b^s_e$ to be finite and invoke the limit $b_e^s \rightarrow \infty$ later. In the reference grid, we have $\hat \theta_e = 0$, so that $\matr M_o = \matr M_c + \Delta \matr M$ with
\begin{align}
    \Delta \matr M &= - b_e^s \hat v_f^2
     ( \vec \mu_e \vec \mu_e^\top + \vec \nu_e \vec \nu_e^\top  )
     \label{eq:DeltaM-switch}
\end{align}
recalling that $\vec \mu_{e} =  \vec \eta_f - \vec \eta_t$ and $\vec \nu_{e} = \vec \zeta_f - \vec \zeta_t$. Using the Woodbury matrix identity we obtain
\begin{align*}
    &\left[ \matr M_c + \Delta \matr M\right]^{-1} 
    =
    \matr M_c^{-1} + 
    \matr M_c^{-1}
    \begin{pmatrix}
        \vec \mu_e & \vec \nu_e
    \end{pmatrix}
    \\
    & \quad \times
    \left[ 
        \frac{1}{b_e^s \hat v_f^2} \eye - 
        \begin{pmatrix}
             \vec \mu_e^\top \\ \vec \nu_e^\top 
        \end{pmatrix}
        \matr M_c^{-1}
        \begin{pmatrix}
            \vec \mu_e & \vec \nu_e
        \end{pmatrix}
    \right]^{-1} 
    \begin{pmatrix}
        \vec \mu_e^\top \\ \vec \nu_e^\top 
    \end{pmatrix}   
    \matr M_c^{-1}.
\end{align*}
At this point we cannot directly impose the limit $b_e^s \rightarrow \infty$ as this implies $ \matr M_c^{-1} \vec \mu_e \rightarrow \vec 0$ and $ \matr M_c^{-1} \vec \nu_e \rightarrow \vec 0$. To circumvent this problem, we rewrite the above expression. Using $\vec \mu_e^\top \vec \nu_e = 0$ and $\vec \mu_e^\top \vec \mu_e =  \vec \nu_e^\top \vec \nu_e = 2$ we find
\begin{align*}
    &  (\eye -  \matr M_c^{-1} \matr M_o) 
    \vec \mu_e
    = - \matr M_c^{-1} \, \Delta \matr M \vec \mu_e
    = - b_{e}^s \hat v_f^2 \matr M_c^{-1} 
    \vec \mu_e
\end{align*}
and an equivalent relation for $\vec \nu_e$. Using these relation as well as their transposes, we obtain
\begin{align}
    & \left[  \frac{1}{b_e^s \, \hat v_f^2} \eye
    -    
    \begin{pmatrix}
        \vec \mu_e^\top \\ \vec \nu_e^\top 
    \end{pmatrix}
    \matr M_c^{-1}
    \begin{pmatrix}
        \vec \mu_e & \vec \nu_e 
    \end{pmatrix}
    \right] 
    \\
    &= 
    \frac{1}{(b_e^s \, \hat v_f^2)^2}
    \begin{pmatrix}
        \vec \mu_e^\top \\ \vec \nu_e^\top 
    \end{pmatrix}
    (\matr M_o - \matr M_o \matr M_c^{-1} \matr M_o) 
    \begin{pmatrix}
        \vec \mu_e & \vec \nu_e 
    \end{pmatrix}.
    \nonumber
\end{align}
Applying the Woodbury formula once more, we finally obtain
\begin{align}
    \matr M_o^{-1} &= 
    \matr M_c^{-1} + 
     \left( \eye - \matr M_c^{-1}  \matr M_o \right)
    \begin{pmatrix}
        \vec \mu_e & \vec \nu_e
    \end{pmatrix} 
    \nonumber \\
    & \quad \times 
    \left[ \begin{pmatrix}
        \vec \mu_e^\top \\ \vec \nu_e^\top 
    \end{pmatrix}
    ( \matr M_o - \matr M_o \matr M_c^{-1} \matr M_o)
    \begin{pmatrix}
        \vec \mu_e & \vec \nu_e
    \end{pmatrix}
    \right]^{-1}
    \nonumber \\
    & \quad \times 
    \begin{pmatrix}
        \vec \mu_e^\top \\ \vec \nu_e^\top 
    \end{pmatrix}
    \left( \eye - \matr M_o \matr M_c^{-1} \right).
    \label{eq:Mo-inverse}
\end{align}
All terms in this expression are finite so that the limit $b_e^s \rightarrow \infty$ becomes trivial. We can thus compute $\matr M_o^{-1}$ from $\matr M_c^{-1}$ by a computationally efficient low-rank matrix update.

We can now compute the state variables in the grid with the open busbar by solving Eq.~\eqref{eq:loadflow-lin-vector}. Substituting Eq.~\eqref{eq:Mo-inverse}, we obtain an expression that is structurally similar to the bus split distribution factors introduced in \cite{van2023bus}. Furthermore, we can carry out an $N-1$ contingency analysis of the open grid topology using results from Sec.~\ref{sec:vlodf}, replacing $\matr M^{-1}$ by $\matr M_o^{-1}$.

\subsection{Example}

We demonstrate the application of the generalized bus split distribution factors. We start from an AC power flow of the base case to obtain the reference values $\hat v_m$ and $\hat \theta_m$ and the matrices $\matr M_c$ and its inverse $\matr M_c^{-1}$. The impact of \emph{all} topology changes can be readily computed via low-rank updates: Finite topology changes can be treated via Eq.~\eqref{eq:woodbury}, the merging of buses and closing of switches via Eq.~\eqref{eq:Minv-merge} and bus splits via Eq.~\eqref{eq:Mo-inverse}. 

\begin{figure}[tb]
    \centering
    \includegraphics[width=0.75\columnwidth]{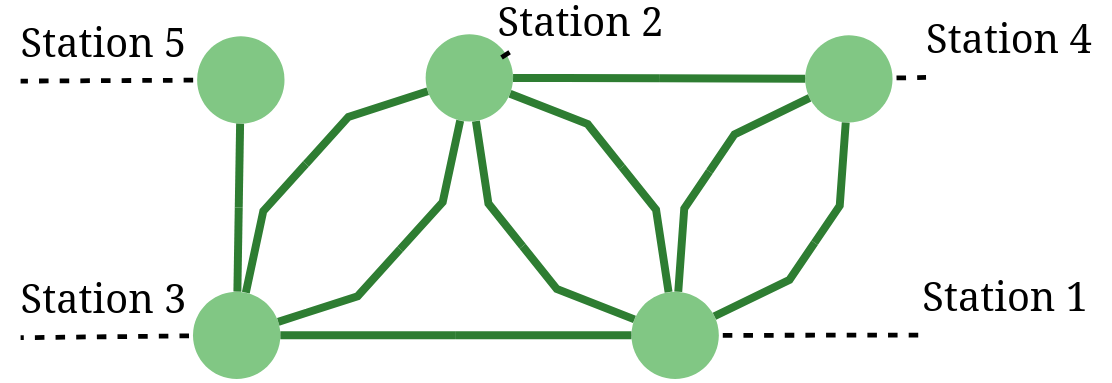}
    \includegraphics[width=0.75\columnwidth]{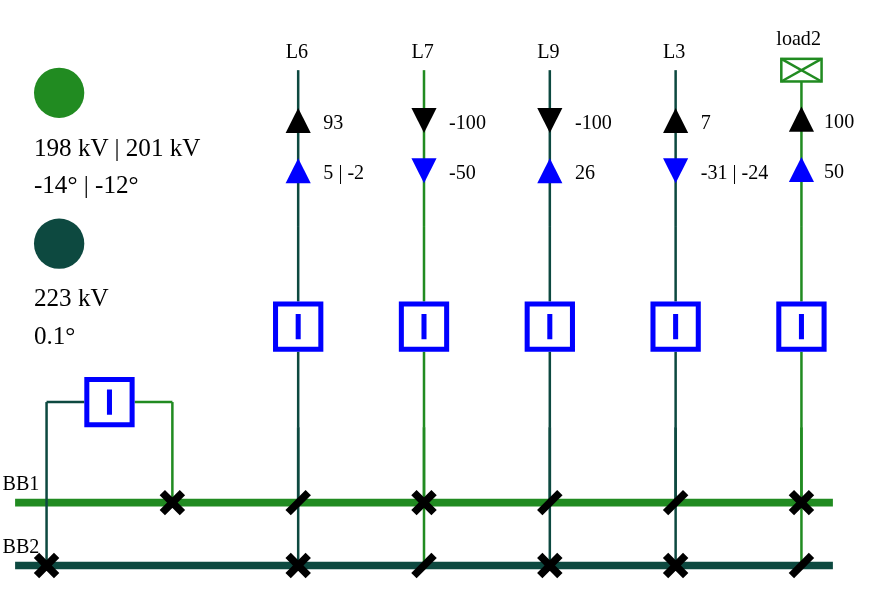}
    \caption{Generalized linear approximation for bus splits. 
    Upper panel: Elementary test grid in the unsplit state. 
    Lower panel: Configuration of the substation 3. Opening the busbar coupler leads to voltage drop and phase angle drop at the busbar 1 (BB1). 
    We provide numerical results for state variables and power flows. Two numbers separated by a $\vert$ compare the generalized DC+ approximation in Eq.~\eqref{eq:Mo-inverse} (left) and full AC results (right).
    AC computations and illustration are produced using PowSyBl~\cite{powsybl2025}. 
    }
    \label{fig:test_grid_bsdf}
\end{figure}

Here, we explore the numerical results and accuracy of the generalized bus split distribution factors for an elementary test (Fig.~\ref{fig:test_grid_bsdf}). We simulate the opening of a busbar coupler at substation VLevel3 comparing the generalized DC+ approximation given by Eq.~\eqref{eq:Mo-inverse} to full AC power flow computation. The split results in a voltage drop at bus bar (BB1) from \SI{225}{kV} to \SI{198}{kV} and a phase angle drop from $\SI{-0.2}{^{\circ}}$ to $\SI{-14}{^{\circ}}$. These values are accurately reproduced by the generalized DC+ approximation that predicts a voltage magnitude drop to \SI{201}{kV} and a voltage angle drop to $\SI{-12}{^{\circ}}$.

\section{Numerical Performance}
\label{sec:performance}

The generalized distribution factors enable new applications in power system operation and planning that require extraordinarily high computational throughput. For example, topology optimization using metaheuristics such as genetic algorithms may involve evaluating hundreds of thousands of switching configurations. When each candidate must additionally undergo a full $N-1$ contingency screening, the resulting workload easily grows to hundreds of millions of power flow evaluations. Similar requirements arise in transmission expansion planning, where many corridor options and operating scenarios must be assessed under large contingency sets, and in operational settings where operators may need to rapidly evaluate large libraries of remedial actions. Importantly, these workloads are inherently parallelizable because each topology and contingency can be evaluated independently, making them well suited for modern GPU architectures.

Meeting these demands is out of reach for nonlinear AC power flow based approaches, whose computational cost is typically orders of magnitude too high for such outer-loop search procedures. The classical DC approximation is computationally much faster and can be implemented very efficiently on GPUs, enabling extremely high evaluation rates~\cite{westerbeck2025accelerated}. However, it cannot capture voltage magnitude variations and reactive power effects and therefore becomes unreliable when voltage security and stability constraints become binding. The generalized distribution factors introduced in this work aim to bridge this gap, combining much of the computational efficiency of linear methods with a significantly improved representation of voltage-sensitive phenomena.

We benchmark the numerical performance of the generalized distribution factors as follows. For various standard test grids we compute the state variables in the $N-0$ case using a standard Newton-Raphson solver. Then we do an $N-1$ contingency analysis using the methods introduced in Sec.~\ref{sec:branch-finite-general}. 
That is, we compute the state variables $\theta_n$ and $v_n$ via Eq.~\eqref{eq:psi-sigma} as well as the currents for all possible non-islanding single branch outages and record the wall time. This is repeated several times to minimize statistical errors. As a figure of merit, we report the number of power flows per second which is defined as 
\begin{align}
    \label{eq:lf-per-sec}
    \mathcal{N} = \frac{n_\mathrm{outages} \times n_\mathrm{repetitions} \times
    \mathrm{batch\_size}
    }{\mathrm{walltime}} \, . 
\end{align}
We set $n_\mathrm{repetitions}=500$ for all GPU experiments and $n_\mathrm{repetitions}=3$ for all CPU experiments.
Parameters of the test grids are summarized in Tab.~\ref{table:benchmark}.

\setlength{\tabcolsep}{5pt}
\begin{table}[t]
\centering
\caption{Benchmark Summary. 
Properties of the test systems and numerical performance for optimized batch size. Performance is computed according to Eq.~\eqref{eq:lf-per-sec} and reported in power flows per second.
}
\label{table:benchmark}
\begin{tabular}{lrrrr}
 Test System & GBnetwork & case2869 & case1354 & case118 \\
$n_{\text{buses}}$        & 2224 & 2869 & 1354 & 118 \\
$n_{\text{outages}}$      & 2083 & 3579 & 1288 & 177 \\
$n_{\text{branches}}$     & 3207 & 4582 & 1991 & 186 \\
\hline
\multicolumn{5}{l}{\textit{Generalized distribution factors (proposed method)}}\\
Load time [s]             & 4.72 & 6.84 & 2.56 & 0.92 \\

Power flows/s & $3.1 \!\times\! 10^6$ & $1.5 \!\times\! 10^6$ & $3.9 \!\times\! 10^6$ & $38.7 \!\times\! 10^6$ \\
Optimal batch size               & 1 & 1 & 128 & 1024 \\
\hline
\multicolumn{5}{l}{\textit{Newton–Raphson (lightsim2grid)}}\\
Power flows/s (1 iter)  & 308 & 192 & 525 & 8688 \\
Power flows/s (20 iter) & 98 & 66 & 176 & 3764 \\
\end{tabular}
\end{table}

We use PyPowSyBl~\cite{powsybl2025} and pandapower~\cite{thurner2018pandapower} for grid file loading and the initial Newton-Raphson process plus some preprocessing routines in native Python. The high-performance $N-1$ loop using the generalized distribution factors is implemented in JAX~\cite{jax}. 
All tests were run on a virtual machine with 48 AMD EPYC 7V13 cores, 440 GB RAM and one NVIDIA A100 GPU with 80GB VRAM. 
Correctness is verified by comparing the results to a single-iteration Newton–Raphson step with the initial guess of the $N-0$ result. 
Our comparison baseline is lightsim2grid \cite{lightsim2grid} on a single CPU core. We run either one or 20 Newton-Raphson iterations with an initial guess from the $N-0$ computation. We emphasize that running one iteration is comparable in accuracy to the generalized distribution factors introduced in this article. 

The numerical results demonstrate that the generalized distribution factors lend themselves well to massively parallel execution on modern GPU architectures. Even for large transmission networks with several thousand buses, our implementation achieves more than one million power flow evaluations per second, while smaller systems reach several tens of millions of evaluations per second. Compared to the AC-based benchmark lightsim2grid, the proposed method improves numerical throughput by several orders of magnitude. When compared to a single Newton–Raphson iteration, which yields an approximation of comparable accuracy, the speedup reaches roughly four orders of magnitude, primarily due to the massive parallelism enabled by GPU execution. Compared to a full Newton–Raphson solve with 20 iterations, corresponding to a fully converged AC power flow solution, the performance gap becomes even larger. 
The approach comes at the cost of several seconds of preprocessing for matrix construction and grid loading before the first results can be generated, making it unsuitable for workloads with fewer than a few hundred cases.
This substantial speedup highlights the advantage of linearized distribution-factor approaches for large-scale contingency screening and topology analysis.

As expected, the computational performance is lower than that of the plain DC formulation due to the additional state variables and the larger Jacobian matrix involved. In earlier work on GPU-accelerated DC distribution factors, we reported throughputs of up to billions of power flows per second under comparable workloads~\cite{westerbeck2025accelerated}. Nevertheless, the present method retains much of the computational efficiency of DC-based techniques while providing a substantially improved representation of voltage-sensitive phenomena.

Overall, these results show that the proposed generalized distribution factor formulation, combined with efficient GPU parallelization, enables extremely fast screening of contingencies and topology modifications. This capability opens the door to new applications in topology optimization and transmission planning, where large numbers of candidate topologies must be evaluated under extensive contingency sets.

\section{Conclusion and Outlook}
\label{sec:conclusion}

We have introduced a framework for voltage-sensitive distribution factors based on a linearized AC power flow formulation. This approach enables efficient evaluation of contingency scenarios and dynamic topology modifications, including line outages and busbar operations. The DC+ framework enables further extensions as the formulation of the matrix $\matr M$ aligns well with fixed Jacobian approaches.

The current implementation demonstrates the feasibility of approximating voltage and reactive power responses with good accuracy, while providing more accurate active power and voltage angle responses than the traditional DC approximation. While the linearization introduces limitations, particularly in reactive power estimation, the results show strong alignment with full AC simulations across different test grids. Further studies suggest that deviations are largest when the voltages violate operational limits. Even in these cases, DC+ provides valuable insights as it typically predicts the occurrence of violations while it cannot capture the precise magnitude.

The DC+ framework is based on low-rank updates via the Woodbury matrix identity, making it well-suited for GPU acceleration~\cite{westerbeck2025accelerated}. 
Hence, it is a promising candidate for time-critical applications such as heuristic algorithms for grid topology optimization.
In our GPU implementation we achieve more than one million power flow evaluations per second even for large transmission networks with several thousand buses. These results demonstrate that voltage-sensitive contingency screening and topology modification analysis can be performed at extremely high throughput, enabling large-scale search and optimization procedures that would be computationally prohibitive with nonlinear AC power flow methods.



\end{document}